\documentclass[12pt]{article}
\usepackage{amsmath}
\usepackage{graphicx}
\usepackage{enumerate}
\usepackage{natbib}
\usepackage{url} 

\newcommand{\blind}{1}

\usepackage{fix-cm}
\usepackage{url}
\usepackage[export]{adjustbox}
\usepackage{amssymb}
\usepackage{amsfonts,bm}
\usepackage{amsmath}
\usepackage[figuresright]{rotating}
\usepackage[toc,page]{appendix}
\usepackage{color}
\usepackage{natbib}
\usepackage{mathrsfs}
\usepackage{booktabs}
\usepackage{enumitem}
\setlist{nolistsep}
\usepackage{array}
\usepackage{graphicx}
\usepackage[ruled,linesnumbered]{algorithm2e}
\usepackage{float}
\usepackage{epstopdf}
\usepackage{mathtools}
\usepackage{subfigure}
\usepackage{setspace}
\usepackage{xr}
\graphicspath{{figure/}}

\newtheorem{remark}{Remark}
\newtheorem{lemma}{Lemma}
\newtheorem{proposition}{Proposition}

\newtheorem{assumption}{Assumption}
\newtheorem{definition}{Definition}

\numberwithin{theorem}{section}
\renewcommand{\baselinestretch}{1.25}
\newcommand{\pr}{^{\prime}}

\newcommand{\n}{^{(n)}}

\newcommand{\bth}{{\bf \theta}}

\newcommand{\thetab}{{\pmb \theta}}

\DeclareMathOperator*{\argmax}{arg\,max}
\DeclareMathOperator*{\argmin}{arg\,min}
\DeclareMathAlphabet\mathbfcal{OMS}{cmsy}{b}{n}

\newcommand{\Zb}{{\bf Z} }

\def\bth{\mbox{\boldmath$\theta$}}

\def\bepsilon{\mbox{\boldmath$\epsilon$}}

\def\r{{\bf r}}

\def\u{{\bf u}}

\def\0{{\bf 0}}

\def\Y{{\bf Y}}

\def\I{{\bf I}}
\def\X{{\bf X}}

\def\0{{\bf 0}}

\def\1{{\bf 1}}

\newcommand{\Yb}{{\bf Y} }
\newcommand{\Xb}{{\bf X} }
\newcommand{\vb}{{\bf v} }

\begin{document}

\def\spacingset#1{\renewcommand{\baselinestretch}%
{#1}\small\normalsize} \spacingset{1}

\if1\blind
{
  \title{\bf Quantiles, Ranks and Signs in Metric Spaces}
  \author{Hang Liu 
  , Xueqin Wang, Jin Zhu, Heping Zhang \\
    University of Science and Technology of China \\
    Sun Yat-Sen University \\
    Yale University \\
    }
    \date{}
  \maketitle
} \fi

\if0\blind
{
  \bigskip
  \bigskip
  \bigskip
  \begin{center}
    {\LARGE\bf  Quantiles, Ranks and Signs in Metric Spaces}
\end{center}
  \medskip
} \fi

\bigskip
\begin{abstract}
  Non-Euclidean data {become more prevalent in practice}, necessitating the development of  {a framework for statistical inference analogous to that for Euclidean data. Quantile is one of the most important concepts in traditional statistical inference; we introduce the counterpart, both locally and globally, for data objects in metric spaces.} This is realized by expanding upon the metric distribution function proposed by \citet{Wang2021}. Rank and sign are defined at local and global levels as a natural consequence of the center-outward ordering of metric spaces brought about by the local and global quantiles. The theoretical properties are established, such as the root-$n$ consistency and uniform consistency of the local and global empirical quantiles and the distribution-freeness of ranks and signs. The empirical metric median, which is defined here as the 0th empirical global metric quantile, is proven to be resistant to contamination by means of both theoretical and numerical approaches. Quantiles have been shown to be valuable through extensive simulations in a number of metric spaces. Moreover, we introduce a family of fast rank-based independence tests for a generic metric space.   Monte Carlo experiments show good finite-sample performance of the test. 
\end{abstract}

\noindent%
{\it Keywords:}  Metric space, Local metric quantile, Global metric quantile, Metric rank, Metric sign
\vfill

\newpage
\spacingset{1.9} 
\section{Introduction}

Data in non-Euclidean spaces can be found in various fields. Examples include directional data \citep{Marinucci2008, Marinucci2011}, symmetric positive definite (SPD) matrices \citep{smith2013functional}, samples of probability distribution functions in Wasserstein spaces \citep{petersen2021wasserstein}, to name only a few. In order to perform nonparametric statistical inference and analysis for such data, we need to develop  {fundamental}  concepts, including quantiles, ranks, and signs for non-Euclidean data.

 Unlike the real line, however, there is no canonical ordering in non-Euclidean spaces, which makes 
  {it is challenging to develop necessary and important concepts, methods, and theories.} 
To tackle this issue, various types of depths have been proposed to introduce ordering in non-Euclidean spaces, including the angular Tukey depth \citep{Liu1992} and distance-based depth \citep{Pandolfo2018} both for spherical data,  the intrinsic zonoid depth for SPD matrices \citep{Chau2019}, the halfspace depth for curve data \citep{de2021depth}, and so on. These notions of depths, however, focus only on a specific metric space. However, relevant literature on general metric spaces is scarce.  By defining halfspace through the underlying distance metric, \cite{Dai2022} proposed a notion of depth at a point as the least probability measure of the halfspaces containing this point. Like the halfspace depth in the Euclidean space, the depth of \cite{Dai2022} is computationally inefficient. Another type of depth is the depth profile proposed by \cite{dubey2022depth}, referring to the distribution of the distances between a point and the other elements in the metric space.

Different from \cite{Dai2022} and \cite{dubey2022depth}, where depths are defined but concepts of quantiles are missing, this paper proposes novel concepts of quantiles from both the local and global perspectives for general metric spaces. Our definition is based on the metric distribution function (MDF), which is the probability measure of a closed ball characterized by two points in metric space \citep{Wang2021}. According to the MDF levels, our local quantiles at a specific location in the metric space are defined.
In order to define quantiles from a global perspective, we take a quantity resulting from aggregating the local information of the MDF for all points in the metric space and define the global quantiles based on the levels of the distribution function of this quantity.  The empirical versions of the local and global quantiles are defined as the discrete sample versions of their population counterparts. Observing the local and global quantiles induces center-outward orderings in the metric space locally and globally, respectively, we therefore define the local and global ranks and signs based on their empirical quantile levels.

In addition to attractive interpretation, as described previously, our novel concepts of the local and global quantiles, ranks, and signs also enjoy a number of desirable theoretical properties. For example, the empirical local and global quantiles are root-$n$ and uniformly consistent; the local and global ranks and signs are distribution-free. Compared with the depths in \cite{Dai2022} and \cite{dubey2022depth},  our empirical global quantiles are easy to implement in practice since the computation is straightforward, does not involve any optimization process, and is computationally efficient with complexity $O(n^2\log(n))$ when the global quantiles are computed for all the $n$ sample points. As illustrated via numerical study in Section~\ref{Sec.Simulation}, distributions in various metric spaces are well-characterized by the global quantiles, and the empirical quantiles have good finite-sample performance, manifesting the reasonableness of our novel concepts. Moreover, theoretical analysis of the breakdown point and simulation study under various metric spaces show that our {\it empirical metric median}, defined as the $0$th empirical global metric quantile, is more robust against contamination than the {\it metric halfspace median} of  \cite{Dai2022} and the {\it transport median} of \cite{dubey2022depth}.\parskip=3pt 

The global ranks and signs in this paper naturally allow us to define rank-based tests for a generic metric space, which have seldom been considered in the literature on rank-based statistical inference. As an illustration, in Section~\ref{Sec.ranktest}, we consider a class of rank-based two-sample independence tests, which are computationally efficient with a complexity of $O(n^2\log(n))$. Monte Carlo experiments show the consistency of our test and good finite-sample performance even under heavy-tailed distributions.

This paper is organized as follows. Section~\ref{sec.MDF} briefly introduces the MDF, its empirical version, and the asymptotic properties in \cite{Wang2021}, which are the building blocks of the main concepts of this paper.  Sections~\ref{sec.local} and \ref{Sec.GlobalQuan}, from the local and global perspectives, respectively,  define the metric quantile, rank, and sign, and establish the asymptotic results, including the root-$n$ and uniform consistency, of the empirical metric quantile. Moreover, computation aspects and robustness of the global metric quantile are discussed in Section~\ref{Sec.GlobalQuan}. Section~\ref{Sec.ranktest} proposes the metric rank-based independence test. Section~\ref{Sec.Simulation} illustrates the concepts of the global metric quantile and its empirical version via Monte Carlo experiments under various metric spaces. Numerical analysis comparing the robustness of the empirical metric median and the counterparts of \cite{Dai2022} and \cite{dubey2022depth} is also included in Section~\ref{Sec.Simulation}.  The finite-sample performance of the metric rank-based independence test is investigated via Monte Carlo experiments in Section~\ref{Sec.Simulation}.
Section~\ref{Sec.Conclusion} concludes the paper and provides some perspectives for future research. All proofs, some preliminaries, simulation details, additional results, and a discussion of the computation algorithm are collected in Appendices.

\section{Preliminaries}\label{sec.MDF}

\citet{Wang2021} proposed the concept of the {\it metric distribution function} (MDF) in a metric space, serving as a quasi-distribution in metric space. In this section, we briefly introduce the MDF and its empirical version, which are important ingredients for defining quantiles in the metric space.

\subsection{Metric distribution function}

An ordered pair $({\cal M}, d)$ is called {\it metric space}  if ${\cal M}$ is a set and $d$ is a metric or distance on ${\cal M}$ (see Appendix~\ref{sec.example} for some examples of non-Euclidean metric spaces that have wide applications in various fields). Throughout the paper, we denote by $\bar{B}(\u, r) := \{\vb: d(\u, \vb) \leq r\}$ the closed ball with the center $\u$ and radius $r \geq 0$ on $({\cal M}, d)$, and by $B(\u, r) := \{\vb: d(\u, \vb) < r\}$ the corresponding open ball. 


Denote by $\mu$ a probability measure associated with a random variable $\X \in {\cal M}$. 
Let
$\delta(\u, \vb, \X) := I(d(\u, \X) \leq d(\u, \vb)),$ 
where $I(\cdot)$ is the indicator function.  The metric distribution function  of $\mu$ on ${\cal M}$  measures the probability that $\X \in {\cal M}$ is contained in $\bar{B}(\u, r)$ and is defined as 
$$F_{\mu}^{\cal M}(\u, \vb) \coloneqq \mu(\bar{B}(\u, \r)) 
= {\rm E}(\delta(\u, \vb, \X))$$
for $\u, \vb \in {\cal M}.$
Different from the traditional distribution function in an Euclidean space, which is a function of one variable and is defined from a global perspective, $F_{\mu}^{\cal M}$ is a probability distribution function defined locally (at $\u \in {\cal M}$) for the closed ball $\bar{B}(\u, r) \subset {\cal M}$ characterized by two variables $\u$ and~$\vb$.

\subsection{Empirical metric distribution function}
Based on the definition of the MDF, its empirical version can naturally be defined as follows.  Let $\X\n \coloneqq  \{\X_1, \ldots, \X_n\}$ be an i.i.d. sample   generated from the probability measure $\mu$ on ${\cal M}$.  For  $\u, \vb \in {\cal M}$, the empirical metric distribution function (EMDF) is the empirical probability measure of the sample being contained in $\bar{B}(\u, r)$. It is defined as 
\begin{equation}\label{Def.EmpDF}
F_{\mu, n}^{\cal M}(\u, \vb) \coloneqq  \frac{1}{n} \sum_{i=1}^n \delta(\u, \vb, \X_i).
\end{equation}

Similar to the empirical distribution function on the real line, $F_{\mu, n}^{\cal M}$ also satisfies the Glivenko-Cantelli and Donsker asymptotic, which are important properties for conducting nonparametric statistical inference in a metric space and for establishing asymptotics of the quantiles proposed in this paper. See Section~\ref{App.GCDonsker} in the Supplementary material for details of the Glivenko-Cantelli and Donsker properties. The establishment of the two properties of the EMDF relies on the following assumptions imposed in \citet{Wang2021}.

\begin{assumption}\label{ass.GC}
$\mu$ on ${\cal M}$ satisfies that
$n^{-1} {\rm E}_{\X}[\log({\rm card}({\cal F}(\X\n)))] \rightarrow 0$, 
where ${\rm card}$ denotes the cardinality of a set and ${\cal F}(\X\n) \coloneqq \{h(\X_1), \ldots, h(\X_n) \vert h \in {\cal F}\}$.
\end{assumption}
%
%
%
\begin{assumption}\label{ass.VC}
The collection of the indicator functions of closed balls on ${\cal M}$, i.e. ${\cal F} \coloneqq \{\delta(\u, \vb, \cdot):  \u, \vb \in {\cal M}\}$, is a VC class\footnote{See \cite{vaart1996weak} for the definition.} with a finite VC-dimension. 
\end{assumption}


%
%

\section{Local metric quantile and its empirical version}\label{sec.local}
\subsection{Local metric quantile}
Once we have the definition of the MDF for the metric space $({\cal M}, d)$, the quantile can naturally be defined according to the level of the MDF.
 Note that, in contrast to the traditional distribution function in an Euclidean space, the MDF is a probability distribution function defined locally at $\u \in {\cal M}$; hence, the quantile in Definition~\ref{Def.Local} is also defined from a local perspective for each $\u \in {\cal M}$. It is a measure of outlyingness with respect to a fixed point in the metric space.\vspace{-3mm}

\begin{definition}\label{Def.Local}
Given a point $\u \in {\cal M}$, the $\tau$th quantile ($\tau \in [0, 1]$) of $\mu$ on ${\cal M}$ is 
$q^{\cal M} (\u, \tau) \coloneqq \{\vb \in {\cal M}: F_{\mu}^{\cal M}(\u, \vb) = \tau\}$,
the $\tau$th {\it quantile region} is $\bar{q}^{\cal M} (\u, \tau) \coloneqq \{\vb \in {\cal M}: F_{\mu}^{\cal M}(\u, \vb) \leq \tau\}$, 
and the $\tau$th {\it inner quantile region} is 
$\tilde{q}^{\cal M} (\u, \tau) \coloneqq \{\vb \in {\cal M}: F_{\mu}^{\cal M}(\u, \vb) < \tau\}$.
\end{definition}

Throughout, we make the following mild assumption on $\mu$, which entails some desirable properties of the quantiles such as nestedness and $\sqrt{n}$-consistency of the empirical local metric quantile defined in Section~\ref{Sec.EmpLocal}. 

\begin{assumption}\label{ass.density}
$\mu$ is absolutely continuous with density function $f$, which is positive and continuous everywhere on ${\cal M}$. 
\end{assumption}

In analog to the quantile regions on the real line, the quantile regions in Definition~\ref{Def.Local} enjoy the following property of nestedness. It implies that the quantile $q^{\cal M} (\u, \tau)$ induces a center-outward ordering of nested regions, with $\u$ being the center. See Section~\ref{sec.proofs} Supplementary material for its proof. \vspace{-2mm}

\begin{proposition}\label{Prop.nest}
If Assumption~\ref{ass.density} holds, then for $0 \leq \tau_1 < \tau_2 \leq 1$ and any $\u \in {\cal M}$, $\bar{q}^{\cal M} (\u, \tau_1) \subsetneq  \bar{q}^{\cal M} (\u, \tau_2)$. \vspace{-2mm}

\end{proposition}

\subsection{Empirical local metric quantile}\label{Sec.EmpLocal}

Given an i.i.d. sample $\X\n  \coloneqq  \{\X_1, \ldots, \X_n\}$  generated from the probability measure $\mu$ on ${\cal M}$, similar to Definition~\ref{Def.Local}, we can naturally define the empirical version of the local metric quantile based on the EMDF. It measures the outlyingness of the sample with respect to a fixed point in the metric space. \vspace{-3mm}

\begin{definition}\label{Def.EmpLocal}
Given a point $\u \in {\cal M}$, the $\tau$th empirical quantile ($\tau \in [0, 1]$) of $\mu$ on ${\cal M}$ is
\begin{equation*}
q_n^{\cal M} (\u, \tau) \coloneqq 
\begin{cases}
\u, & \tau \in [0, \frac{1}{n}) \\
\argmin_{\X_i  \in \X\n}  F_{\mu, n}^{\cal M}(\u, \X_i) \quad {\rm s.t.} \quad F_{\mu, n}^{\cal M}(\u, \X_i)  \geq \tau, & \tau \in [\frac{1}{n}, 1].
\end{cases}
\end{equation*}
\end{definition}

\begin{remark}
Note that the $\tau$th population quantile in Definition~\ref{Def.Local} yields a set of contours in ${\cal M}$, whereas the empirical version, with probability one, is a point in ${\cal M}$  due to the absolute continuity of $\mu$. The $\tau$th sample quantile may alternatively be defined as a set of points through factorizing $n$, which would yield a set of contours in ${\cal M}$ after performing interpolation. Specifically, by factorization of $n = n_1 n_2 + n_0$, where $n_1, n_2 \rightarrow \infty$ as $n \rightarrow \infty$ and $n_0$, as the multiplicity of the $0$th sample quantile, is a small number, one can define the  $\tau$th ($\tau \geq 1/n_2$) sample quantile as a set of $n_1$ points according to the ordering of  $F_{\mu, n}^{\cal M}(\u, \X_i), i = 1, \ldots, n$. This approach is equivalent to defining $\tau$th sample quantile as the set of $[\tau, \tau + n_1/n)$th sample quantiles in Definition~\ref{Def.EmpLocal}. This type of definition, however, usually requires a relatively large $n$ (especially when the dimension is large) to obtain a regular contour. In the context of Euclidean space, defining the sample contour through the factorization of $n$ can be found in \cite{Hallinetal2021}. A similar remark also applies to the global empirical metric quantile defined in Section~\ref{Sec.EmpGlobal}.  \vspace{-3mm}


\end{remark}

Proposition~\ref{Prop.ConsistLocal} states that the $\tau$th empirical local metric quantile, with probability one, lies in a small neighborhood (of order $n^{-1/2}$) of the $\tau$th local quantile contours, which is generated by the difference of two nested inner quantile regions. See Section~\ref{sec.proofs} in the Supplementary material for its proof.

\begin{proposition}\label{Prop.ConsistLocal}
For any $\u \in {\cal M}$, when $\tau \in [0, \frac{1}{n})$, $q_n^{\cal M} (\u, \tau) = q^{\cal M}(\u, 0)$; when $\tau \in [\frac{1}{n}, 1]$, under Assumptions~\ref{ass.VC} and~\ref{ass.density}, with probability one, \vspace{-3mm}
$$q_n^{\cal M} (\u, \tau) \subset \tilde{q}^{\cal M} (\u, \tau + \frac{1}{n} + \frac{1}{\sqrt{n}}   \mathbb{G}_{\mu, n}^{\cal M}(\u, q_n^{\cal M} (\u, \tau))) \setminus \tilde{q}^{\cal M} (\u, \tau + \frac{1}{\sqrt{n}} \mathbb{G}_{\mu, n}^{\cal M}(\u, q_n^{\cal M} (\u, \tau))),$$ \vspace{-3mm}
where $\mathbb{G}_{\mu, n}^{\cal M}(\u, \vb) := \sqrt{n}(F_{\mu, n}^{\cal M}(\u, \vb) - F_{\mu}^{\cal M}(\u, \vb)), \, \u, \vb \in {\cal M}$.

\end{proposition}

The following result, moreover, states the uniform consistency (in terms of $\u$ and $\tau$) of the empirical local metric quantile. See Section~\ref{sec.proofs} in the Supplementary material for its proof.

\begin{proposition}\label{Prop.unif.local}
If Assumption~\ref{ass.GC} holds, then
$\underset{n \rightarrow \infty}{\lim} \underset{\u \in {\cal M}, \tau \in [0, 1]}{\sup} {\rm P}(q_n^{\cal M} (\u, \tau) \not\subset q^{\cal M} (\u, \tau)) = 0.$

\end{proposition}

\subsection{Local metric ranks and signs}\label{sec.local.rank}

On a real line, the distribution function, as a monotonic increasing map from $\mathbb{R}$ to $[0, 1]$, introduces the left-to-right order relation. Its empirical version, also as a monotonic increasing map from the i.i.d. sample to the regular grid on $[0, 1]$, naturally allows us to define ranks of a sample by multiplication of the sample size $n$ (or $n+1$). In our setting of a general metric space, when a point $\u \in {\cal M}$ is given, we can similarly define ranks locally by the local empirical distribution function. Moreover, signs can also be defined locally according to the sign of $F_{\mu, n}^{\cal M}(\u, \X_j) n/(n+1) - 1/2, 1\leq j \leq n.$  \vspace{-3mm}


\begin{definition}
Let $\{\X_1, \ldots, \X_n\}$ be an i.i.d. sample  generated from the probability measure $\mu$ on ${\cal M}$. Then given a point $\u \in {\cal M}$, the local metric rank of $\X_j, 1\leq j \leq n$ is 
$$R_{\mu, n}^{\cal M; {\rm loc}}(\u, \X_j) := n F_{\mu, n}^{\cal M}(\u, \X_j) = \sum_{i=1}^n \delta(\u, \X_j, \X_i)$$
and the local metric sign of $\X_j$ is
$$S_{\mu, n}^{\cal M; {\rm loc}}(\u, \X_j) := {\rm sign}\left(\frac{n F_{\mu, n}^{\cal M}(\u, \X_j)}{n+1} - \frac{1}{2}\right) = {\rm sign}\left(\frac{R_{\mu, n}^{\cal M; {\rm loc}}(\u, \X_j)}{n+1} - \frac{1}{2}\right).$$
\end{definition}

An important property that makes ranks and signs successful tools for statistical inference is their distribution-freeness. Similar to their counterparts on the real line, this property also holds for our local metric ranks and signs. See Section~\ref{sec.proofs} in the Supplementary material for its proof. \vspace{-3mm}

\begin{proposition}\label{Prop.DistFree}
Let $\{\X_1, \ldots, \X_n\}$ be an i.i.d. sample  generated from the probability measure $\mu$ on ${\cal M}$. Then given a point $\u \in {\cal M}$,
\begin{enumerate}
\item[(i)] $F_{\mu, n}^{\cal M}(\u, \X\n) := (F_{\mu, n}^{\cal M}(\u, \X_1), F_{\mu, n}^{\cal M}(\u, \X_2), \ldots, F_{\mu, n}^{\cal M}(\u, \X_n))$ is uniformly distributed over the $n!$ permutations of $(1/n, 2/n, \ldots, 1)$;
\item[(ii)] $R_{\mu, n}^{\cal M; {\rm loc}}(\u, \X\n) := (R_{\mu, n}^{\cal M; {\rm loc}}(\u, \X_1), \ldots, R_{\mu, n}^{\cal M; {\rm loc}}(\u, \X_n))$ is uniformly distributed over the $n!$ permutations of $(1, \ldots, n)$;
\item[(iii)] ${\rm P}(S_{\mu, n}^{\cal M; {\rm loc}}(\u, \X_j) = 1) = \frac{\left\lfloor n/2\right\rfloor}{n} = {\rm P}(S_{\mu, n}^{\cal M; {\rm loc}}(\u, \X_j) = -1)$ and ${\rm P}(S_{\mu, n}^{\cal M; {\rm loc}}(\u, \X_j) = 0)   = 1- \frac{2\left\lfloor n/2\right\rfloor}{n}$.
\end{enumerate}

\end{proposition}

\subsection{Isometry-invariance of the local concepts}\label{Sec.EquiLocal}

Now, we discuss the equivariance properties of the distribution functions, local quantiles, ranks, and signs. Let $({\cal M}, d)$ and $({\cal M}\pr, d\pr)$ be two metric spaces, and 
$T: {\cal M} \rightarrow {\cal M}\pr$ be an isometry; that is, 
$d\pr(T\u , T\vb) = d(\u, \vb), \forall \u, \vb \in {\cal M}.$
Proposition~\ref{Prop.eqv.local} states that the distribution function, empirical version, local ranks, and signs are invariant under an isometric transformation. The local quantile and its empirical version are also invariant to an isometric transformation. See Section~\ref{sec.proofs} in the Supplementary material for its proof. \vspace{-5mm}

\begin{proposition}\label{Prop.eqv.local}
Let $\{\X_1, \ldots, \X_n\}$ be i.i.d. samples  generated from the probability measure $\mu$ on ${\cal M}$. Then
\begin{enumerate}
\item[(i)] $F_{\mu}^{{\cal M}\pr}(T\u, T\vb) = F_{\mu}^{\cal M}(\u, \vb)$ and $F_{\mu, n}^{{\cal M}\pr}(T\u, T\vb) = F_{\mu, n}^{\cal M}(\u, \vb)$
\item[(ii)] $q^{{\cal M}\pr} (T\u, \tau) = Tq^{\cal M} (\u, \tau)$ and $q_n^{\cal M\pr} (T\u, \tau) = Tq_n^{\cal M} (\u, \tau)$
\item[(iii)] $R_{\mu, n}^{{\cal M}\pr; {\rm loc}}(T\u, T\X_j) = R_{\mu, n}^{{\cal M}; {\rm loc}}(\u, \X_j)$ and $S_{\mu, n}^{{\cal M}\pr; {\rm loc}}(T\u, T\X_j) = S_{\mu, n}^{{\cal M}; {\rm loc}}(\u, \X_j)$ for $j = 1, \ldots, n$.\vspace{-3mm}
\end{enumerate}

\end{proposition}

\section{Global metric quantile and its empirical version}\label{Sec.GlobalQuan}

\subsection{Global metric quantile}

\begin{figure}[htbp]
\centering
    \includegraphics[scale=0.5]{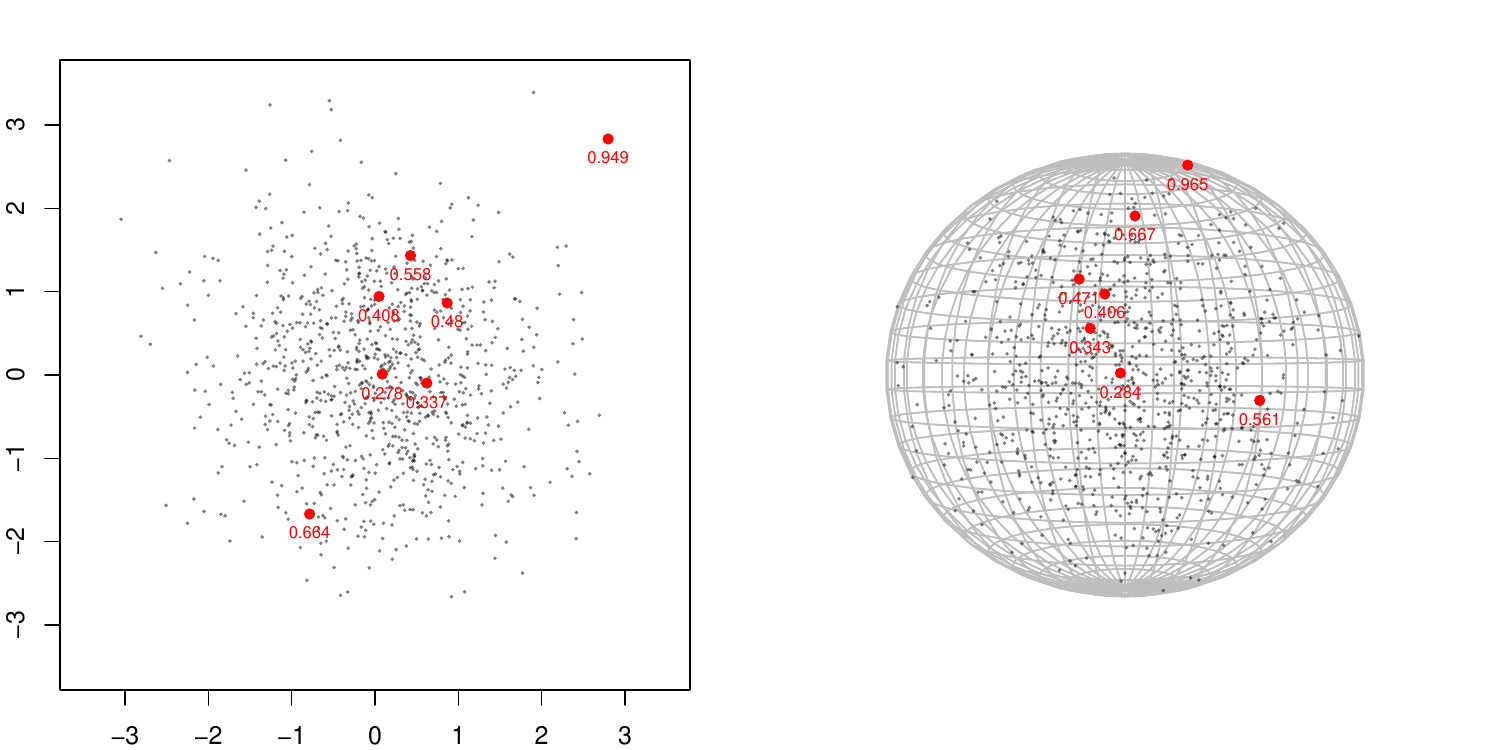}
    \caption{Plots of samples generated from the standard spherical Gaussian distribution in $\mathbb{R}^2$ (left panel) and the von Mises-Fisher distribution in $\mathcal{S}^2$ (right panel), with estimates of $J_{\mu}^{\cal M}(\u)$ marked in red at some points.}\label{Fig.JnR2S2}
\end{figure}

In contrast to the local metric quantile in Section~\ref{sec.local},  this section defines the quantile from a global perspective. Noticing that  $F_{\mu}^{\cal M}(\X, \u)$ provides probability information for $\u$  locally at $\X$,  to define the global quantile, we consider aggregating the local information of $F_{\mu}^{\cal M}(\X, \u)$ for all points $\X \in {\cal M}$. To this end, define $J_{\mu}^{\cal M}(\u) \coloneqq {\rm E}(F_{\mu}^{\cal M}(\X, \u))$ which, 
intuitively, provides an averaged ``population ranking'' information of $\u$ when evaluating at each point in the metric space; in the terminology of sociology, it measures how ``popular'' $\u$ is when voted by all $\X \in {\cal M}$. As a simple illustration of the intuition, we simulate i.i.d. samples from the standard spherical Gaussian in $\mathbb{R}^2$ and a von Mises-Fisher (vMF) distribution in $\mathcal{S}^2$  (see Section~\ref{App.SimSphere} in the Supplementary material for details of the vMF distribution). We compute $J_{\mu, n}^{\cal M}(\u) := n^{-1} \sum_{i=1}^n F_{\mu, n}^{\cal M}(\X_i, \u)$, with $n=1000$, as an estimate of $J_{\mu}^{\cal M}(\u)$.   In Figure~\ref{Fig.JnR2S2} (the left/right panel is for $\mathbb{R}^2$/$\mathcal{S}^2$), values of $J_{\mu, n}^{\cal M}(\u)$ of some points are marked in red. Notice that, in Figure~\ref{Fig.JnR2S2},  $J_{\mu, n}^{\cal M}(\u)$ increases as $\u$ moves away from the center. This motivates our definition of global metric quantile, which essentially relies on the ordering of $J_{\mu}^{\cal M}(\u)$ for each $\u \in {\cal M}$. Letting $L \coloneqq \underset{\u \in {\cal M}}{\min}J_{\mu}^{\cal M}(\u)$ and $U \coloneqq \underset{\u \in {\cal M}}{\max}J_{\mu}^{\cal M}(\u)$, we make the following assumption.

\vspace{-3mm}

\begin{assumption}\label{ass.UL}
$L \neq U$. \vspace{-3mm}
\end{assumption}

We remark that in the case that $L = U$, the global quantiles may simply be defined as the local quantiles $q^{\cal M} (\u, \tau)$ by choosing any point $\u$ on ${\cal M}$, since $J_{\mu}^{\cal M}(\u)$ is uniformly distributed over ${\cal M}$ and any $\u$ would provide the same global information. For instance, when ${\cal M}$ is the unit sphere and $\Xb$ is uniformly distributed, one may choose any $\u \in {\cal M}$ as the $0$th global quantile and define the $\tau$th global quantile as $q^{\cal M} (\u, \tau)$.

Since $\mu$ is absolutely continuous with positive density $f$ everywhere in ${\cal M}$, $J_{\mu}^{\cal M}(\u)$ has a positive density $f_J$ everywhere on $[L, U]$. Denoting by $F_J$ the distribution function of $J_{\mu}^{\cal M}(\u) \in [L, U]$, then $F_J$ is a strict monotonic increasing function of $J_{\mu}^{\cal M}$. We define the global quantiles as follows. \vspace{-3mm}

\begin{definition}\label{Def.global}
The $\tau$th global quantile ($\tau \in [0, 1]$) of $\mu$ on ${\cal M}$ is
$Q^{\cal M} (\tau) \coloneqq \{\u \in {\cal M}: F_J(J_{\mu}^{\cal M}(\u)) = \tau\}$,
the $\tau$th {\it quantile region} is $\bar{Q}^{\cal M} (\tau) \coloneqq \{\u \in {\cal M}: F_J(J_{\mu}^{\cal M}(\u)) \leq \tau\}$, 
and the $\tau$th {\it inner quantile region} is 
$\tilde{Q}^{\cal M} (\tau) \coloneqq \{\u \in {\cal M}: F_J(J_{\mu}^{\cal M}(\u)) < \tau\}.$
\end{definition}

Similar to the local quantile regions, the global quantile regions are also nested, and their intersection is $0$th global quantile region, which is a compact set of $\u \in {\cal M}$ that minimizes $J_{\mu}^{\cal M}(\u)$. This minimizer is called {\it metric median}, which may be unique. For example, when ${\cal M} = \mathcal{S}^p$ and $\mu$ is rotationally symmetric with respect to the pole, the metric median is a parallel or a collection of parallels. Proposition~\ref{Prop.nest.global} states the nestedness property of the global quantile regions. Therefore,  the global quantiles induce a center-outward ordering for the metric space ${\cal M}$. See Section~\ref{sec.proofs} in the Supplementary material for its proof. \vspace{-4mm}

\begin{proposition}\label{Prop.nest.global}
If Assumptions~\ref{ass.density} and~\ref{ass.UL} hold, then for $0 \leq \tau_1 < \tau_2 \leq 1$, $\bar{Q}^{\cal M} (\tau_1) \subsetneq  \bar{Q}^{\cal M} (\tau_2)$. \vspace{-3mm}

\end{proposition}

\subsection{Empirical global metric quantile}\label{Sec.EmpGlobal}

Let $J_{\mu, n}^{\cal M}(\u) = n^{-1} \sum_{i=1}^n F_{\mu, n}^{\cal M}(\X_i, \u)$, which can be considered as an approximation to  $J_{\mu}^{\cal M}(\u) \coloneqq {\rm E}(F_{\mu}^{\cal M}(\X, \u))$, and denote by
$F_{J_n}\n(J_{\mu, n}^{\cal M}(\u)) \coloneqq \frac{1}{n}  \sum_{i=1}^n I(J_{\mu, n}^{\cal M} (\X_i) \leq J_{\mu, n}^{\cal M} (\u))$
the empirical distribution function of $J_{\mu, n}^{\cal M}(\u)$. The definition of empirical global metric quantiles relies on the ordering of $J_{\mu, n}^{\cal M}(\X_i), i = 1, \ldots, n$.

\begin{definition}\label{Def.empGlobal}
The $\tau$th empirical global quantile ($\tau \in [0, 1]$) of $\mu$ on ${\cal M}$ is
\begin{equation*}
Q_n^{\cal M}(\tau) \coloneqq \argmin_{\X_i \in \X\n} F_{J_n}\n(J_{\mu, n}^{\cal M}(\X_i)) \quad {\rm s.t.} \quad F_{J_n}\n(J_{\mu, n}^{\cal M}(\X_i)) \geq \tau.
\end{equation*}

\end{definition}

Proposition~\ref{Prop.ConsistGlobal} states that the $\tau$th empirical global quantile, similar to its local version, is uniformly consistent over $\tau \in [0, 1]$. See Section~\ref{sec.proofs} in the Supplementary material for its proof.

\begin{proposition}\label{Prop.ConsistGlobal}
If Assumptions~\ref{ass.GC}, \ref{ass.density} and~\ref{ass.UL} holds, then
$\underset{n \rightarrow \infty}{\lim} \underset{\tau \in [0, 1]}{\sup} {\rm P}(Q_n^{\cal M} (\tau) \not\subset Q^{\cal M} (\tau)) = 0.$
\end{proposition}

Before introducing the next property of the empirical global quantile, we first state some preliminary results of $J_{\mu}^{\cal M}$ and $J_{\mu, n}^{\cal M}$. Let
$H_{n1}(y) \coloneqq \sqrt{n}\left({\rm P}(J_{\mu, n}^{\cal M}(\X_1) \leq y) - \frac{1}{n} \sum_{i=1}^n I(J_{\mu, n}^{\cal M}(\X_i) \leq y) \right)$
and
$H_{n2}(\u) \coloneqq \sqrt{n}(J_{\mu}^{\cal M}(\u) - J_{\mu, n}^{\cal M}(\u)).$ 
Then we have the following results for $H_{n1}(y)$ and $H_{n2}(\u)$. See Section~\ref{sec.proofs} in the Supplementary material for its proof. \vspace{-3mm}

\begin{lemma}\label{lem.Hn1}
If Assumption~\ref{ass.GC} holds, 
then, as $n\rightarrow \infty$,
$H_{n1}(y) \rightarrow \mathcal{N}(0, F_J(y) - (F_J(y))^2). \vspace{-3mm}$
\end{lemma}
\begin{lemma}\label{lem.Hn2}
If Assumption~\ref{ass.VC} holds, then as $n\rightarrow \infty$,
$H_{n2}(\u) \rightarrow \mathcal{N}(0, {\rm Var}(F_{\mu}^{\cal M}(\X, \u))).$
\end{lemma}

Let 
\begin{equation}\label{eq.Dn}
D_n(\tau) \coloneqq  H_{n1}(J_{\mu, n}^{\cal M}(Q_n^{\cal M} (\tau))) + f_J(J_{\mu}^{\cal M}(Q_n^{\cal M} (\tau)) \left(H_{n2}(Q_n^{\cal M} (\tau)) -  H_{n2}(\X_1)\right).
\end{equation}
Proposition~\ref{Prop.Global} states that the $\tau$th empirical global quantile, similar to its local version, also lies in a small neighborhood (with a radius of order $n^{-1/2}$) of the $\tau$th global quantile contours. See Section~\ref{sec.proofs} in the Supplementary material for its proof. \vspace{-4mm}

\begin{proposition}\label{Prop.Global}
Let Assumptions~\ref{ass.GC}-\ref{ass.UL} hold. Then for any $\tau \in [0, 1]$, with probability one,
\begin{equation}\label{eq.QnQtilde}
Q_n^{\cal M} (\tau) \subset \tilde{Q}^{\cal M} (\tau\n + \frac{1}{n}) \setminus \tilde{Q}^{\cal M} (\tau\n),
\end{equation}
where $\tau\n = \tau + \frac{1}{\sqrt{n}}D_{n}(\tau) + o_{\rm P}(n^{-1/2}),$
and $D_{n}(\tau)$ converges to a normal distribution with mean zero.

\end{proposition}

\subsection{Global metric rank and sign}\label{Sec.ranksign}

The canonical ordering of the real line makes it easy to define ranks and signs in terms of the sample order. In $\mathbb{R}^p, p \geq 2$, however, such ordering does not exist.  Using the measure transportation theory, \cite{Hallinetal2021}  defined ranks and signs based on a center-outward ordering of the sample. In our setting of a general metric space $({\cal M}, d)$, we have defined ranks and signs from a local perspective in Section~\ref{sec.local.rank}. From a global perspective,  since the global quantile regions relying on the ordering of $J_{\mu}^{\cal M} (\u)$ induce a center-outward ordering for ${\cal M}$, ranks and signs will naturally be defined based on the ordering of the empirical version $J_{\mu, n}^{\cal M} (\X_j), 1\leq j \leq n,$ of $J_{\mu}^{\cal M} (\u)$. \vspace{-3mm}

\begin{definition}
The metric rank of $\X_j, 1\leq j \leq n$ is 
$$R_{\mu, n}^{\cal M; {\rm Global}}(\X_j) := n F_{J_n}\n(J_{\mu, n}^{\cal M}(\X_j)) = \sum_{i=1}^n I(J_{\mu, n}^{\cal M} (\X_i) \leq J_{\mu, n}^{\cal M} (\X_j)),$$
and the metric sign of $\X_j$ is
$$S_{\mu, n}^{\cal M; {\rm Global}}(\X_j) := {\rm sign}\left(\frac{n F_{J_n}\n(J_{\mu, n}^{\cal M}(\X_j))}{n+1}  - \frac{1}{2}\right) = {\rm sign}\left(\frac{R_{\mu, n}^{\cal M; {\rm Global}}}{n+1}  - \frac{1}{2}\right).$$
\end{definition}

Proposition~\ref{Prop.DistFreeGlb} states the distribution-freeness property of the global metric rank and sign as an extension of the corresponding property in the Euclidean space $(\mathbb{R}^p, \Vert \cdot \Vert)$ (see Proposition~2.5 in \cite{Hallinetal2021}). See Section~\ref{sec.proofs} in the Supplementary material for its proof.

\begin{proposition}\label{Prop.DistFreeGlb}
Let $\{\X_1, \ldots, \X_n\}$ be an i.i.d. sample  generated from the probability measure $\mu$ on ${\cal M}$. Then
\begin{enumerate}
\item[(i)] $F_{J_n}\n(J_{\mu, n}^{\cal M}(\X\n)) := (F_{J_n}\n(J_{\mu, n}^{\cal M}(\X_1)), F_{J_n}\n(J_{\mu, n}^{\cal M}(\X_2)), \ldots, F_{J_n}\n(J_{\mu, n}^{\cal M}(\X_n)))$ is uniformly distributed over the $n!$ permutations of $(1/n, 2/n, \ldots, 1)$;
\item[(ii)] $R_{\mu, n}^{\cal M; {\rm Global}}(\X\n) := (R_{\mu, n}^{\cal M; {\rm Global}}(\X_1), \ldots, R_{\mu, n}^{\cal M; {\rm Global}}(\X_n))$ is uniformly distributed over the $n!$ permutations of $(1, \ldots, n)$;
\item[(iii)] ${\rm P}(S_{\mu, n}^{\cal M; {\rm Global}}(\X_j) = 1) = \frac{\left\lfloor n/2\right\rfloor}{n} = {\rm P}(S_{\mu, n}^{\cal M; {\rm Global}}(\X_j) = -1)$ and ${\rm P}(S_{\mu, n}^{\cal M; {\rm Global}}(\X_j) = 0)   = 1- \frac{2\left\lfloor n/2\right\rfloor}{n}$.
\end{enumerate}

\end{proposition}

\subsection{Isometry-invariance of the global concepts}\label{Sec.EquiGlobal}

Consider the isometry $T: {\cal M} \rightarrow {\cal M}\pr$ defined in Section~\ref{Sec.EquiLocal}. 
Like their local counterparts, the global quantile, ranks, and signs enjoy the following equivariance properties.  See Section~\ref{sec.proofs} in the Supplementary material for its proof. \vspace{-3mm}

\begin{proposition}\label{Prop.eqv.global}
Let $\{\X_1, \ldots, \X_n\}$ be an i.i.d. sample  generated from the probability measure $\mu$ on ${\cal M}$. Then
\begin{enumerate}
\item[(i)] $J_{\mu}^{{\cal M}\pr}(T\u) = J_{\mu}^{{\cal M}}(\u)$ and $J_{\mu, n}^{{\cal M}}(T\u) = J_{\mu, n}^{{\cal M}\pr}(\u)$;
\item[(ii)] $Q^{{\cal M}\pr} (\tau) = TQ^{\cal M} (\tau)$ and $Q_n^{\cal M\pr} (\tau) = TQ_n^{\cal M} (\tau)$;
\item[(iii)] $R_{\mu, n}^{{\cal M}\pr; {\rm Global}}(T\X_j) = R_{\mu, n}^{{\cal M}; {\rm Global}}(\X_j)$ and $S_{\mu, n}^{{\cal M}\pr; {\rm Global}}(T\X_j) = S_{\mu, n}^{{\cal M}; {\rm Global}}(\X_j)$ for $j = 1, \ldots, n$.
\end{enumerate}

\end{proposition}

\subsection{Relation with depth function}

Various types of depth functions have been proposed in order to tackle the problem of non-canonical ordering in multi-dimensional Euclidean spaces and non-Euclidean spaces, including halfspace depth \citep{Tukey1975}, simplicial depth \citep{Liu1990}, spatial depth \citep{serfling2002depth}, to name just a few,  in  Euclidean space, and halfspace depth \citep{Dai2022} and depth profiles \citep{dubey2022depth} in non-Euclidean spaces. In the context of Euclidean spaces, \cite{Zuo2000} formally defined the statistical depth function which should possess four properties: {\it (i)} affine invariance with respect to some linear transformation; {\it (ii)}  the depth function is maximized at the center; {\it (iii)} the depth function decreases monotonically along any fixed ray through the center; {\it (iv)} the depth function approaches zero as the point approaches infinity. 

As mentioned in Section~\ref{Sec.GlobalQuan}, our global quantiles hinge on the ordering of $J_{\mu}^{\cal M}(\u)$. Let $\bth \in {\cal M}$ denote a point in the metric median $Q^{\cal M}(0)$ of $({\cal M}, d)$. Proposition~\ref{Prop.depth} states that under certain conditions of the metric space and the probability measure, $D_{\mu}^{\cal M}(\u):= 1 - F_J(J_{\mu}^{\cal M}(\u))$ can be interpreted as a ``metric depth'' function satisfying the above properties of a depth function in an  Euclidean space. See Section~\ref{sec.proofs} in the Supplementary material for its proof. \vspace{-3mm}

\begin{proposition}\label{Prop.depth}
The metric depth function $D_{\mu}^{\cal M}(\u)$ satisfies
\begin{enumerate}
\item[(i)] $D_{\mu}^{{\cal M}\pr}(T\u) = D_{\mu}^{\cal M}(\u)$ under an isometric transformation $T: {\cal M} \rightarrow {\cal M}\pr$ defined in Section~\ref{Sec.EquiGlobal};
\item[(ii)] $D_{\mu}^{\cal M}(\u)$ attains the maximum value at any $\bth \in Q^{\cal M}(0)$;
\item[(iii)] when $({\cal M}, d)$ is a geodesic space and the probability measure $\mu$ satisfies $J_{\mu}^{\cal M}(\gamma(t_1)) \leq J_{\mu}^{\cal M}(\gamma(t_2))$ for any $0\leq t_1 \leq t_2 \leq 1$, where  $\gamma: [0, 1] \rightarrow  {\cal M}$ is an arbitrary geodesic curve starting from $\bth$, i.e. $\gamma(0) = \bth \in Q^{\cal M}(0)$, then $D_{\mu}^{\cal M}(\gamma(t_2)) \leq D_{\mu}^{\cal M}(\gamma(t_1))$;
\item[(iv)] $D_{\mu}^{\cal M}(\u) = 0$ for any $\u \in  \underset{\vb \in {\cal M}}{\argmax}J_{\mu}^{\cal M}(\vb)$.
\end{enumerate}

\end{proposition}

Suppose $({\cal M}, d)$ is the geodesic space as in part (iii) of Proposition~\ref{Prop.depth}. Let 
$$\Gamma_{\gamma(t_1), \gamma(t_2))} := \{\X \in {\cal M}: d(\X, \gamma(t_2)) \geq d(\X, \gamma(t_1))\} \vspace{-2mm}$$
and  
$$\widetilde{\Gamma}_{\gamma(t_1), \gamma(t_2))} := \{\X \in {\cal M}: d(\X, \gamma(t_2)) < d(\X, \gamma(t_1))\} \vspace{-2mm}$$
be two sets that partition $({\cal M}, d)$ according to the distance to $\gamma(t_1)$ and $\gamma(t_2)$. Then, the condition in (iii) is equivalent to 
$$\int_{\Gamma_{\gamma(t_1), \gamma(t_2))}} \mu(B(\X, \gamma(t_2) \setminus B(\X, \gamma(t_1)) {\rm d}\mu \geq \int_{\widetilde{\Gamma}_{\gamma(t_1), \gamma(t_2))}} \mu(B(\X, \gamma(t_1) \setminus B(\X, \gamma(t_2)) {\rm d}\mu$$
for any geodesic curve $\gamma \in {\cal M}$ and any $0\leq t_1 \leq t_2 \leq 1$.  Loosely speaking, this means that points that are close to the metric center along a geodesic curve have large densities. This condition holds in Euclidean space for, e.g., elliptical distributions with the density function decreasing monotonically along any fixed ray through the mean. \vspace{-4mm}


\subsection{Computation aspects}\label{Sec:computation}

We flag that an advantage of our empirical global metric quantiles (EGMQ) is that they are easy to implement in computation and computationally efficient, with a time complexity of $O(n^2\log{n})$.\footnote{We highlight that the EGMQ is more computationally efficient than the depth of \cite{Dai2022}, of which the time complexity is $O(n^3)$.} In Section~\ref{App.algorithm} in the Supplementary material, we provide a more detailed discussion of the computation aspects and algorithm.

\subsection{Breakdown point of the empirical metric median}

The breakdown point of an estimator is the minimal proportion of contamination in a sample that can make this estimator arbitrarily bad, and it is a tool quantifying the robustness of this estimator.  In this section, we study the breakdown point of the empirical metric median (EMM), defined as the $0$th empirical global metric quantile, i.e., 
$$\thetab_n^{\Xb} \coloneqq Q_n^{\cal M}(0)  = \argmin_{\Xb_i \in \Xb\n} J_{\mu, n}^{\cal M}(\X_i).$$ 
In particular, suppose the sample 
$\Zb^{(m+n)} \coloneqq  \{\Xb_1, \ldots, \Xb_n, \Yb_1, \ldots, \Yb_m\}$ contains $m$ contaminated observations $\Yb^{(m)} \coloneqq  \{\Yb_1, \ldots, \Yb_m\}$. Let $\Zb_i, i = 1, \ldots, m+n$ be an element of $\Zb^{(m+n)}$ and denote by $\tilde{F}_{\mu, m+n}^\mathcal{M}(\u, \vb)$ the EMDF of $\Zb^{(m+n)}$. Let 
$\thetab_{m, n}^{\Xb, \Yb} \coloneqq \argmin_{\Zb_i \in \Zb^{(m+n)}} \tilde{J}_{\mu, m+n}^{\cal M}(\Zb_i),$
where $\tilde{J}_{\mu, m+n}^{\cal M}(\u)\coloneqq (m+n)^{-1} \sum_{i=1}^{m+n} \tilde{F}_{\mu, m+n}^\mathcal{M}(\Zb_i, \u)$. Then the {\it finite sample breakdown point} of the empirical metric median is defined as
\begin{equation}\label{def.breakdown}
\epsilon^*_n \coloneqq \min_m \left\lbrace \frac{m}{m+n}: \sup_{\Yb^{(m)}} d(\thetab_n^{\Xb},  \thetab_{m, n}^{\Xb, \Yb}) = \infty \right\rbrace,
\end{equation}
with the convention that $\epsilon^*_n = 1$ if the RHS of \eqref{def.breakdown} is an empty set (for instance, when ${\cal M}$ is bounded). The following proposition gives a lower bound of $\epsilon^*_n$. See Section~\ref{sec.proofs} in the Supplementary material for the proof. \vspace{-3mm}

\begin{proposition}\label{Prop.breakdown}
The finite sample breakdown point of the empirical metric median satisfies
\begin{equation}\label{eq.fsbp}
\epsilon^*_n \geq \frac{1 - n^{-1}\sum_{i=1}^{n} F_{\mu, n}^\mathcal{M}(\Xb_i, \thetab_n^{\Xb})}{2 - n^{-1}\sum_{i=1}^{n} F_{\mu, n}^\mathcal{M}(\Xb_i, \thetab_n^{\Xb})}.
\end{equation}
Consequently, when Assumption~\ref{ass.GC}, \ref{ass.density} and~\ref{ass.UL} holds, as $n\rightarrow \infty$, the {\it asymptotic breakdown point} $\epsilon^* := \underset{n\rightarrow \infty}{\lim}\epsilon^*_n$ satisfies that with probability one,
\begin{equation}\label{eq.asbp}
\epsilon^* \geq 
\frac{1 - {\rm E}(F_{\mu}^\mathcal{M}(\Xb, \thetab^{\Xb}))}{2 -  {\rm E}(F_{\mu}^\mathcal{M}(\Xb, \thetab^{\Xb}))},
\end{equation}
where $\thetab^{\Xb}  \coloneqq Q^{\cal M}(0)  = \argmin_{\Xb \in {\cal M}} J_{\mu}^{\cal M}(\Xb)$ is the {\it metric median} of ${\cal M}$.
\end{proposition}
\vspace{-4mm}

\begin{remark}
Proposition~\ref{Prop.breakdown} implies that the lower bound of $\epsilon^*$ depends on the value of $J_{\mu}^{\cal M}(\thetab^{\Xb}) = {\rm E}(F_{\mu}^\mathcal{M}(\Xb, \thetab^{\Xb}))$: a distribution with a smaller value of $J_{\mu}^{\cal M}(\thetab^{\Xb})$ has a larger value of the lower bound. For any symmetric distribution, one can show that ${\rm E}(F_{\mu}^\mathcal{M}(\Xb, \thetab^{\Xb})) < 1/2$ when Assumption~\ref{ass.density} holds, and hence this lower bound is always greater than $1/3$, which is the lower bound of the {\it metric halfspace median} (MHM) of \cite{Dai2022}. Table 1 lists the lower bounds of the asymptotic breakdown points \footnote{The asymptotic values are obtained from the empirical counterparts with $n=1000$.}of the EMM and MHM under the spherical Gaussian distribution and skew-$t_6$ distribution in $\mathbb{R}^2$, the von Mises-Fisher (vMF) and tangent vMF distribution in the unit sphere $\mathcal{S}^2$, the Wishart distribution in SPD(3) with degrees of freedom $3$ and scaling matrix
$\begin{pmatrix}
    1 & 0.6 & 0.36 \\
    0.6 & 1 & 0.6 \\
    0.36 & 0.6 & 1
\end{pmatrix};$
see Section~\ref{Sec.Simulation} and Section~\ref{sec:SimDetails} in the Supplementary material for more details of these distributions. Notice that the lower bound of the breakdown point of the EMM is uniformly greater than that of the MHM. For the space of SPD matrices, the difference between the EMM and MHM is more significant, with the value of the EMM close to $0.5$ and the value of the MHM being $0.2$. This phenomenon of greater robustness of the EMM than the MHM for the space of SPD matrices is further confirmed via simulation study in Section~\ref{sec.simRobust}.
\end{remark}

\begin{table} 
\caption{The lower bounds of the asymptotic breakdown points of the EMM and MHM under various distributions.}\label{Tab.breakdown}
    \centering
    \begin{tabular}{c|ccccc}
        \toprule
            & Gaussian & skew-$t_6$ & vMF  & tangent vMF & Wishart \\
        \midrule
        EMM & 0.42     & 0.42       & 0.41 & 0.42        & 0.47    \\
        MHM & 0.33     & 0.32       & 0.33 & 0.32        & 0.21    \\
        \bottomrule
    \end{tabular}
\end{table}

\section{Metric rank-based independence test}\label{Sec.ranktest}

Thanks to the distribution-free property, rank-based tests have been applied quite successfully in Euclidean spaces for various statistical models, including location \citep{hodges1956efficiency},   regression models \citep{jureckova1971nonparametric, koul1971asymptotic, hallin2022efficient}, autoregressive time series \citep{koul1993r},  
non-linear time series \citep{mukherjee2007, liu2022r}, to name a few. In a generic metric space, due to the lack of the concept of rank, their application has been much less considered. Based on the global metric rank defined in Section~\ref{Sec.ranksign}, we are able to construct rank-based tests for a generic metric space. In this section, we illustrate this point by considering a class of rank-based tests of independence.

\subsection{Test statistics}

Let $\X \in ({\cal M}_1, d_1)$ and $\Y \in ({\cal M}_2, d_2)$ be two random variables with distributions satisfying Assumption~\ref{ass.density}. The problem of interest is to test the null
$\mathcal{H}_0: \X \,\, \text{and} \,\, \Y \, \text{are independent}$ 
based on $n$ independent copies $(\X_1, \Y_1), \ldots, (\X_n, \Y_n)$ of $(\X, \Y)$. Let $R_{ni}^{\X}$ and $R_{ni}^{\Y}$ denote the global metric ranks of $\X_i$ and $\Y_i$ among $\X_1, \ldots, \X_n$ and $\Y_1, \ldots, \Y_n$, respectively. The rank-based test of independence we consider is based on the classical linear rank statistic (see, e.g., \citet[Chapter~4]{Hajek}), which takes the form
$$T\n_{\varphi_1, \varphi_2} := \sum_{i=1}^n \varphi_1\left(\frac{R_{ni}^{\X}}{n+1}\right) \varphi_2\left(\frac{R_{ni}^{\Y}}{n+1}\right),$$
where $\varphi_1: [0, 1) \rightarrow \mathbb{R}$ and $\varphi_2: [0, 1) \rightarrow \mathbb{R}$ are {\it score functions} satisfying the following traditional assumption. \vspace{-3mm}

\begin{assumption}\label{Ass.scores}
The score functions $\varphi_1$ and $\varphi_2$ {\it (i)} are square-integrable, that is, $\int_0^1 \varphi_{\ell}(u) {\rm d}u <\infty, \ell = 1, 2$, and  {\it (ii)} are differences of  two continuous monotonic increasing functions.
\end{assumption}
\vspace{-4mm}

Assumption~\ref{Ass.scores} is quite mild and is satisfied by all square-integrable functions with bounded variation. The metric rank-based test statistic takes the form
$$W\n_{\varphi_1, \varphi_2} := (T\n_{\varphi_1, \varphi_2} - {\rm E}(T\n_{\varphi_1, \varphi_2}))/({\rm Var}(T\n_{\varphi_1, \varphi_2}))^{1/2},$$
where the expectation and variance are taken under $\mathcal{H}_0$. 
 In view of Theorem~6.1.8.1 in \cite{Hajek}, $W\n_{\varphi_1, \varphi_2}$ is asymptotically standard normal under $\mathcal{H}_0$ as $n\rightarrow \infty$. 

The test statistic $W\n_{\varphi_1, \varphi_2}$ relies on the choice of score functions $\varphi_1$ and $\varphi_2$. Below, we provide an example of a standard score function that is widely employed in the Euclidean space (see, e.g., \cite{Hajek} for the univariate case and \cite{HLL2020} for the multivariate case). 

\textbf{Example} ({\it Spearman} score).   A simple choice is $\varphi_1(u) = u = \varphi_2(u)$. The corresponding linear rank statistic is 
$T\n_{\text{Sp}} = \sum_{i=1}^n R_{ni}^{\X} R_{ni}^{\Y}/(n+1)^2$
and the test statistic \vspace{-3mm}
$$W\n_{\text{Sp}} = \frac{12}{n(n+1)\sqrt{n-1}}\sum_{i=1}^n R_{ni}^{\X} R_{ni}^{\Y} -  \frac{3(n+1)}{\sqrt{n-1}} \vspace{-3mm}$$
upto a scale factor $(n-1)^{-1/2}$, coincides with the the Spearman rank correlation.


\section{Numerical study}\label{Sec.Simulation}

\subsection{Numerical illustration of the global metric quantile}\label{Sec.simQuantile}

In this section, we illustrate our concepts of the global quantile and its empirical version through Monte Carlo experiments on the following metric spaces; for the sake of conciseness of this section, more details of the simulation for these metric spaces are deferred to Section~\ref{sec:SimDetails}  in the Supplementary material. For each scenario, we generate an i.i.d. sample of size $n$ and plot its empirical global metric quantiles. The purpose of these experiments is to show the reasonableness of the concepts under various settings and to explore the finite-sample performance of the empirical quantiles. 
In order to plot $\tau$th population quantile, which is a set of contours, we simulate an i.i.d. sample of a large size $N$, 
and we use \textsf{geom\_contour\_filled} function in the \textsf{ggplot2} package \citep{wickham2016ggplot2} to estimate and plot the contour.
We flag the small-sample global quantiles in red triangles in some plots to show the consistency of the sample global quantiles.

\begin{figure}[h]
     \centering
     \begin{subfigure} 
         \centering
         \includegraphics[scale=0.4]{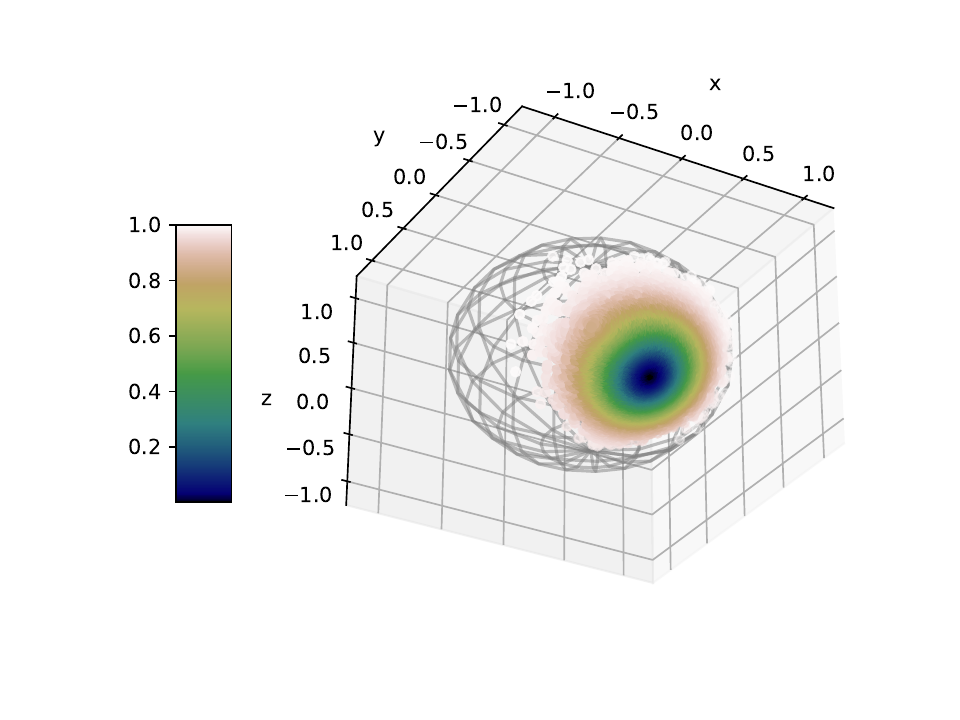}
     \end{subfigure}
     \begin{subfigure} 
         \centering
         \includegraphics[scale=0.4]{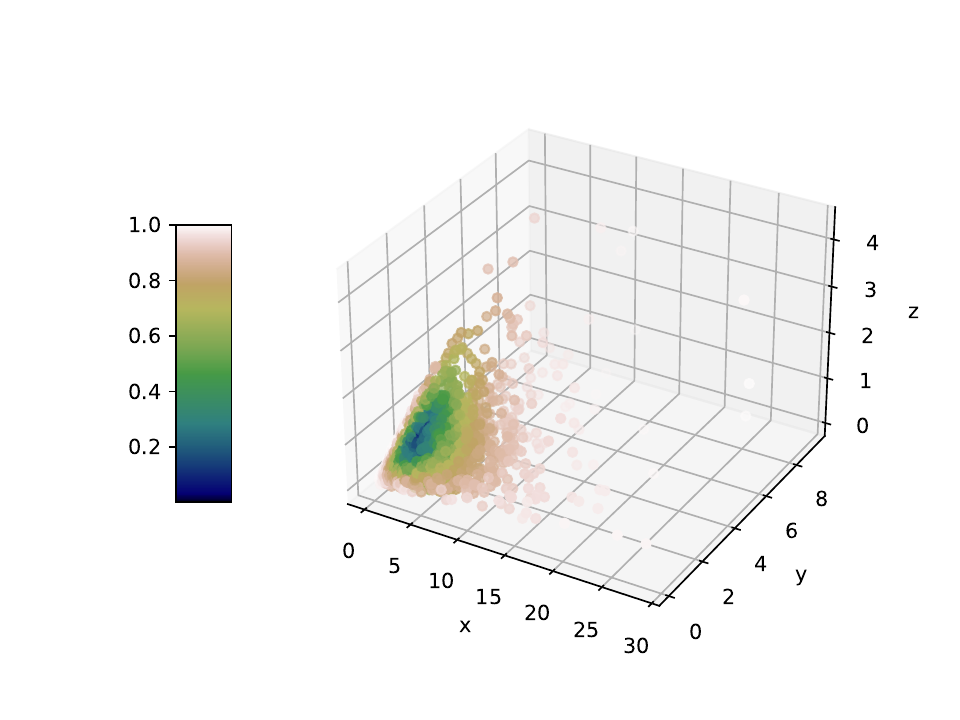} \vspace{-3mm}
     \end{subfigure}
        \caption{Left panel:  the tangent von Mises-Fisher distribution in the unit sphere $\mathcal{S}^2$. Right panel: space of SPD matrices where each component of the matrices is generated from a log-normal distribution (each point in the plot represents the lower diagonal values of an SPD matrix).
The color reflects the level of the global metric quantile.}
       \label{Fig:S2tanvMF}
\end{figure}

\begin{enumerate}
\item[(1)] The Euclidean space $\mathbb{R}^2$, in which we consider two types of distributions as follows, and we set $n = 500$, $N = 10000$. \color{black}
\begin{itemize}
\item[(i)] The spherical Gaussian distribution $\mathcal{N}(\0, \I_2)$ (see Figure~\ref{R2Gauss} in Section~\ref{App.simQuantile}  in the Supplementary material).
\item[(ii)] Skew-$t_6$ distribution (see Figure~\ref{R2SkewT}  in Section~\ref{App.simQuantile}  in the Supplementary material).
\end{itemize}

\item[(2)] The unit sphere $\mathcal{S}^2$, where three types of distributions on $\mathcal{S}^2$ are considered. We set $n = 500$, $N = 10000$. \color{black}
\begin{itemize}
\item[(i)] vMF distribution (see Figure~\ref{Fig:S2vMF}  in Section~\ref{App.simQuantile} in the Supplementary material). 
\item[(ii)] tangent vMF distribution. (see the left panel of Figure~\ref{Fig:S2tanvMF}).
\item[(iii)] mixture of two vMFs (see Figure~\ref{Fig:S2Mixture}  in Section~\ref{App.simQuantile}  in the Supplementary material). 
\end{itemize}

\item[(3)] The Wasserstein space $W_2$ equipped with the Wasserstein distance defined in equation~\eqref{eq:WassDist} with $p=2$. Here, each sample point $\Xb_i$ is a Gaussian distribution $\mathcal{N}(0, \sigma_i^2)$ with $\sigma_i$ being generated from the following distributions. We set $n = 200$, $N = 3000$. \color{black} Scatterplots of $\sigma_i$ for $i = 1, \ldots, N$, the corresponding quantile contours, and the small-sample quantiles are demonstrated in Figures~\ref{Fig:WassUnif} and \ref{Fig:WassBeta}  in Section~\ref{App.simQuantile}  in the Supplementary material for the cases (i) and (ii), respectively. 
\begin{itemize}
\item[(i)] The uniform distribution $U(0, 1)$. 
\item[(ii)] The beta distribution Beta$(2, 5)$.
\end{itemize}


\item[(4)] Space of symmetric positive definite (SPD) matrices equipped with the affine-invariant Riemannian metric defined in equation~\eqref{eq:SPDdist}. Letting $\X_i = ((x_{i}, y_{i})^\top, (y_i, z_i)^\top)$, $x_{i}, y_{i}, z_{i}$ are generated from a log-normal distribution. We set $n = 500$, $N = 3000$. \color{black} Scatterplots of $(x_i, y_i, z_i)^\top$ and the corresponding metric quantiles are demonstrated in the right panel of Figure~\ref{Fig:S2tanvMF}. 

\item[(5)] The collection of phylogenetic trees equipped with the BHV distance. 
We consider the simplest tree space $\mathbb{T}^{3}$ with 
three leaves and one interior edge. In our simulation, 80\% trees lie on the second branch, and 10\% trees lie on the remaining branch, respectively. We choose $n = 500$, $N = 12000$. 
Figure~\ref{fig:tree_rbeta}  in Section~\ref{App.simQuantile}  in the Supplementary material presents scatterplots for the random tree observations and the corresponding metric quantiles. \vspace{-2mm}

\end{enumerate}

For the Euclidean space, we observe that the population quantile contours nicely conform to the shape of the underlying distribution for the spherical Gaussian and skew-$t_3$ cases, where both the skewness and heavy tails are presented. Moreover, the empirical global quantiles demonstrate excellent finite-sample performance as they all lie close to their population quantiles, which is consistent with Proposition~\ref{Prop.ConsistGlobal} and  Proposition~\ref{Prop.Global}.

A similar phenomenon is observed in Figures~\ref{Fig:S2vMF}-\ref{Fig:S2Mixture} for the unit sphere. 
For the Wasserstein space, when the sample is uniform over $(W_2, d)$ (Figure~\ref{Fig:WassUnif}), the global metric quantile contours also are nearly center-outward uniformly distributed; when the sample are skew distributed (Figure~\ref{Fig:WassBeta}), the quantile contours also present a similar pattern of skewness. For the space of SPD matrices, the right panel of Figure~\ref{Fig:S2tanvMF} shows a center-outward ordering shape of the metric quantiles.  Regarding the space of phylogenetic trees, we observe that all the low-level quantiles are located on the second branch, and this is consistent with our simulation setting where  80\% trees lie on the second branch.

\subsection{Numerical evidence for robustness}\label{sec.simRobust}

In this section, we investigate the robustness of the empirical metric median (EMM) through Monte Carlo experiments. 
The metric halfspace median (MHM) of \cite{Dai2022}, which was shown to be robust against contamination, and the estimated {\it transport median} (ETM) of \cite{dubey2022depth}
will be used for comparison. We consider the Euclidean space, spherical space, and the space of SPD matrices. For each space, the i.i.d. sample of size $n=100$ is generated from a ``contaminated" distribution $(1 - \alpha) {\rm P}_1 + \alpha {\rm P}_2$, 
where $\alpha \in [0, 0.5)$ denotes the proportion of the contamination, and ${\rm P}_1$ and ${\rm P}_2$ are as follows.

\begin{enumerate}
    \item[(i)] The Euclidean space $\mathbb{R}^3$.  ${\rm P}_1$ is the standard spherical Gaussian distribution $\mathcal{N}(\mathbf{0}, \mathbf{I}_3)$, and  
    ${\rm P}_2$ is the multivariate Cauchy distribution in which each dimension is sampled from a Cauchy distribution with zero location and unit scale; 
    \item[(ii)] The spherical space $\mathcal{S}^2$. ${\rm P}_1$ and ${\rm P}_2$ are the vMF distributions, both with concentration parameter 1, but with different locations $(1, 0, 0)^\top$ and $(1/\sqrt{3}, 1/\sqrt{3}, 1/\sqrt{3})^\top$;
    \item[(iii)] The space of SPD matrices SPD($3$). ${\rm P}_1$ and ${\rm P}_2$ are the Wishart distributions, both with degrees of freedom 3, but with different scaling matrices $\mathbf{I}_{3}$ and $\begin{pmatrix}
        1 & 0.6 & 0.36 \\
        0.6 & 1 & 0.6 \\
        0.36 & 0.6 & 1
    \end{pmatrix}.$
\end{enumerate}

In each experiment, we compute the distances of the EMM, MHM, and ETM to the center of ${\rm P}_1$. The experiment is repeated 100 times for each metric space, and we compare the robustness of the EMM, MHM, and ETM through their averaged distances to the center of ${\rm P}_1$. Plots of the averaged distances for different levels of $\alpha$ are demonstrated in Figures~\ref{fig:euc_mean_distance}, \ref{fig:sphere_mean_distance}  and~\ref{fig:spd_mean_distance} in Section~\ref{App.simRobust}  in the Supplementary material for the simulation schemes (i), (ii) and (iii), respectively. For the compact metric space $\mathcal{S}^2$, the EMM, MHM, and ETM show similar performance. When it comes to the unbounded spaces $\mathbb{R}^3$ and SPD($3$), their differences become more distinct. In particular, for $\mathbb{R}^3$, the MHM is the least robust one, and the EMM and ETM demonstrate similar robustness.  For the SPD($3$), the EMM is the most robust one, and as $\alpha$ gets close to $0.5$, the gaps between their averaged distances become much more significant. This is also inconsistent with our analysis of the breakdown points of the EMM and MHM in Table~\ref{Tab.breakdown}.\vspace{-3mm}

\subsection{Numerical study of the metric rank-based test}\label{Sec.simRankTest}

In this section, we investigate via Monte Carlo experiments the finite-sample performance of the metric rank-based independence test. We compare our test with the {\it ball covariance}-based (BCov) independence test \citep{pan2019ball}, and the {\it distance covariance}-based (DCov) independence test \citep{lyons2013distance}. 
We simulate $n$ independent copies $\{(\X_i, \Y_i)\}_{i=1}^n$ by letting $$\Y_i = (k\X_i + 0.8\bepsilon_i + 0.5\I_2) (k\X_i + 0.8\bepsilon_i + 0.5\I_2)\pr, i = 1, \ldots, n,$$ 
where $k$ is a positive value, $\X_i$ is generated from the bivariate spherical Gaussian distribution ${\cal N}(\0, \I_2)$, and $\bepsilon_i$, independent of $\X_i$, is generated from ${\cal N}(\0, \I_2)$ or the bivariate Cauchy distribution. Hence, $\X$ is in the Euclidean space, and $\Y$ is in the space of SPD matrices.

\begin{figure}[htbp]
    \centering
    \begin{subfigure} 
        \centering
        \includegraphics[width=0.47\textwidth]{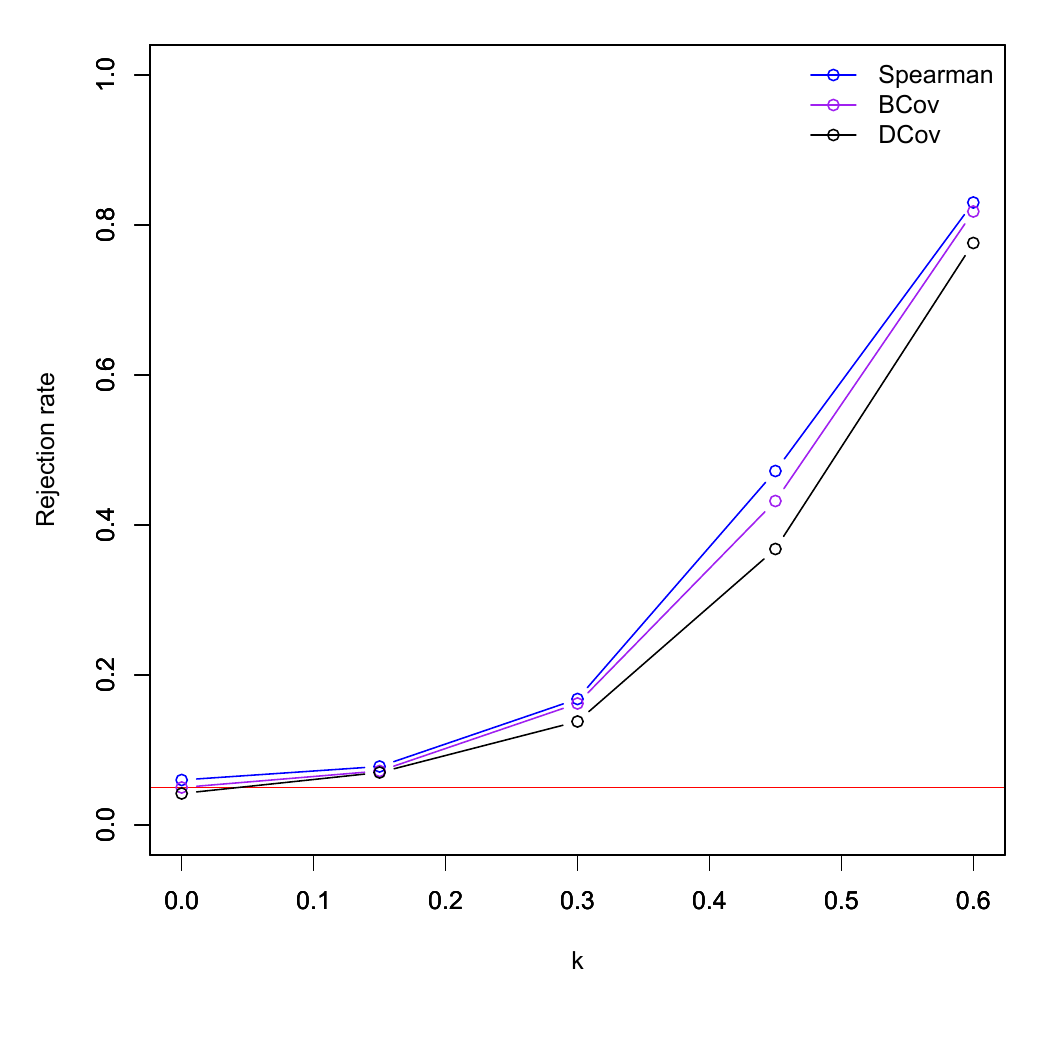}
    \end{subfigure}
    \begin{subfigure} 
        \centering
        \includegraphics[width=0.47\textwidth]{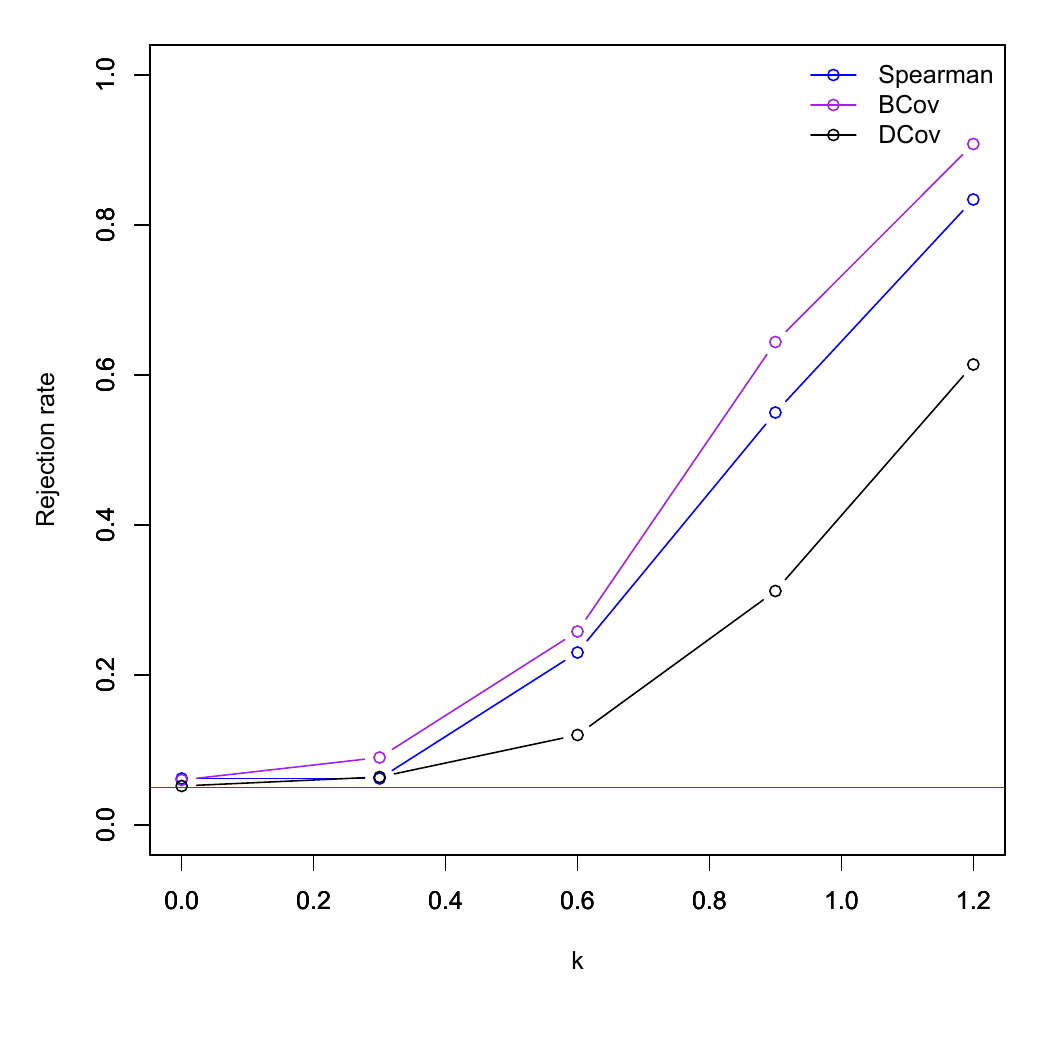}
    \end{subfigure}
    \vspace*{-30pt}
    \caption{The rejection rates of the Spearman metric rank-based, BCov and DCov tests for the spherical normal distribution (left panel) and Cauchy distribution (right panel) of $\bepsilon_i$ for different values of $k$. The red line represents the nominal level.}
  \label{plot.kDiff}
\end{figure}

The rejection rates of the Spearman metric rank-based, BCov and DCov tests at $5\%$ nominal level for different values of $k$ are shown in Figure~\ref{plot.kDiff}, where the left-panel corresponds to the case that $\bepsilon_i$ is generated from the Gaussian distribution and the right-panel is for the Cauchy distribution. For both distributions, all the tests have rejection rates close to the nominal level when $\X$ and $\Y$ are independent (e.g., $k=0$). All the tests have similar power under the Gaussian distribution. Under the Cauchy distribution, however, the Spearman and BCov tests have similar rejection rates, which are significantly greater than that of the DCov test. This illustrates the Spearman and BCov tests (the BCov test is robust since it is constructed from the metric distribution function). 

Next, we compare the rejection rates of these tests by fixing $k$ and changing the sample size $n$. For the Gaussian distribution of $\bepsilon_i$, we set $k=0.8$, and the results are shown in the left panel of Figure~\ref{plot.nDiff}  in Section~\ref{App.simRankTest}  in the Supplementary material; for the Cauchy distribution, we set $k=2$, and the results are given in the right panel of Figure~\ref{plot.nDiff}. For both distributions, rejection rates of all the tests get close to one as $n$ increases. As in Figure~\ref{plot.kDiff}, the Spearman and BCov tests have similar performance, and both have rejection rates remarkably greater than DCov.

In addition to the great finite-sample performance, our metric rank-based test is also computationally more efficient than the BCov and DCov tests. Specifically, our test has the complexity of  $O(n^2\log(n))$. In contrast, the BCov and DCov tests use $R$ times permutations replications to approximate the distribution under the null hypothesis, having complexities $O(R n^2\log(n))$ and $O(R n^2)$, respectively. Thus, the proposed test is more computationally attractive than the benchmarked tests for metric-valued data.

\begin{figure}[htbp]
    \centering
    \begin{subfigure} 
        \centering
        \includegraphics[width=0.47\textwidth]{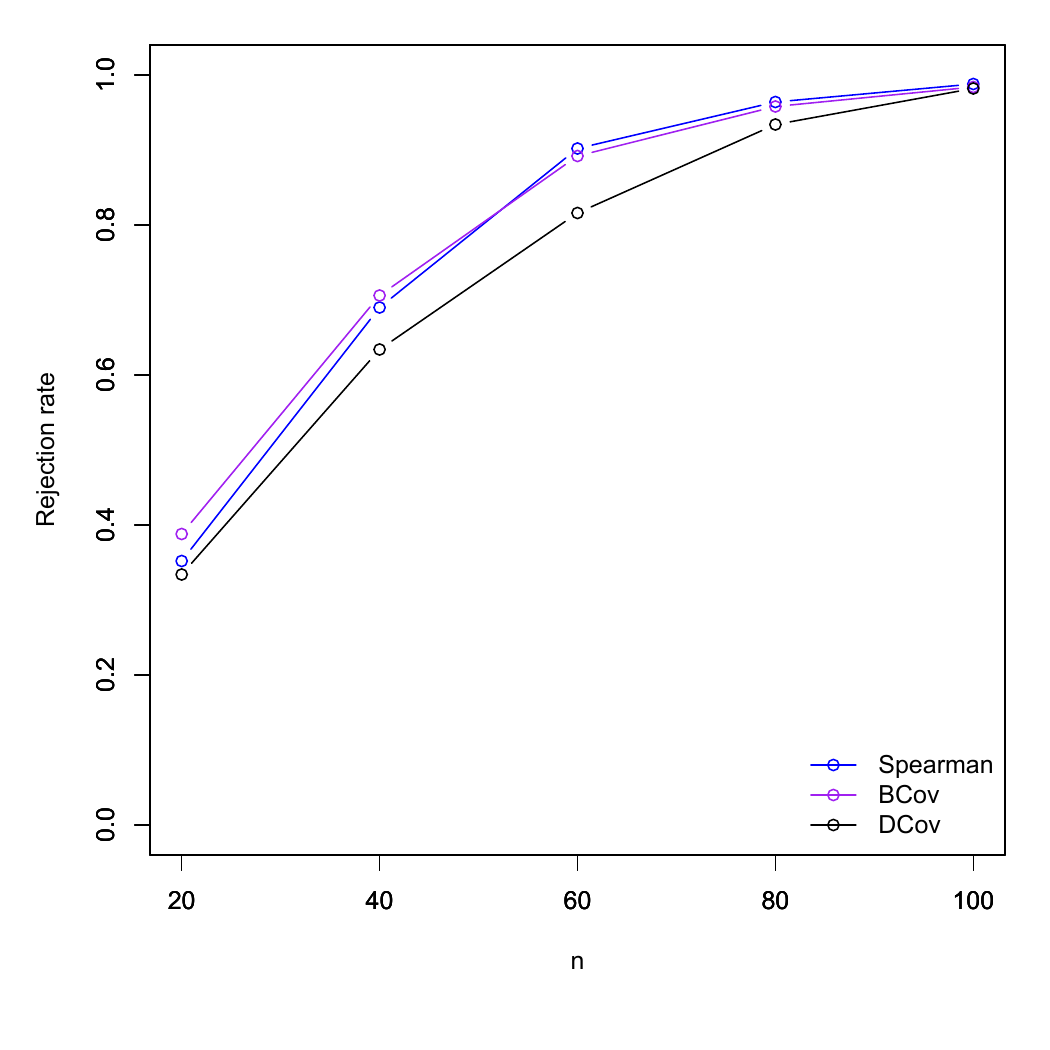}
    \end{subfigure}
    \begin{subfigure}
        \centering
        \includegraphics[width=0.47\textwidth]{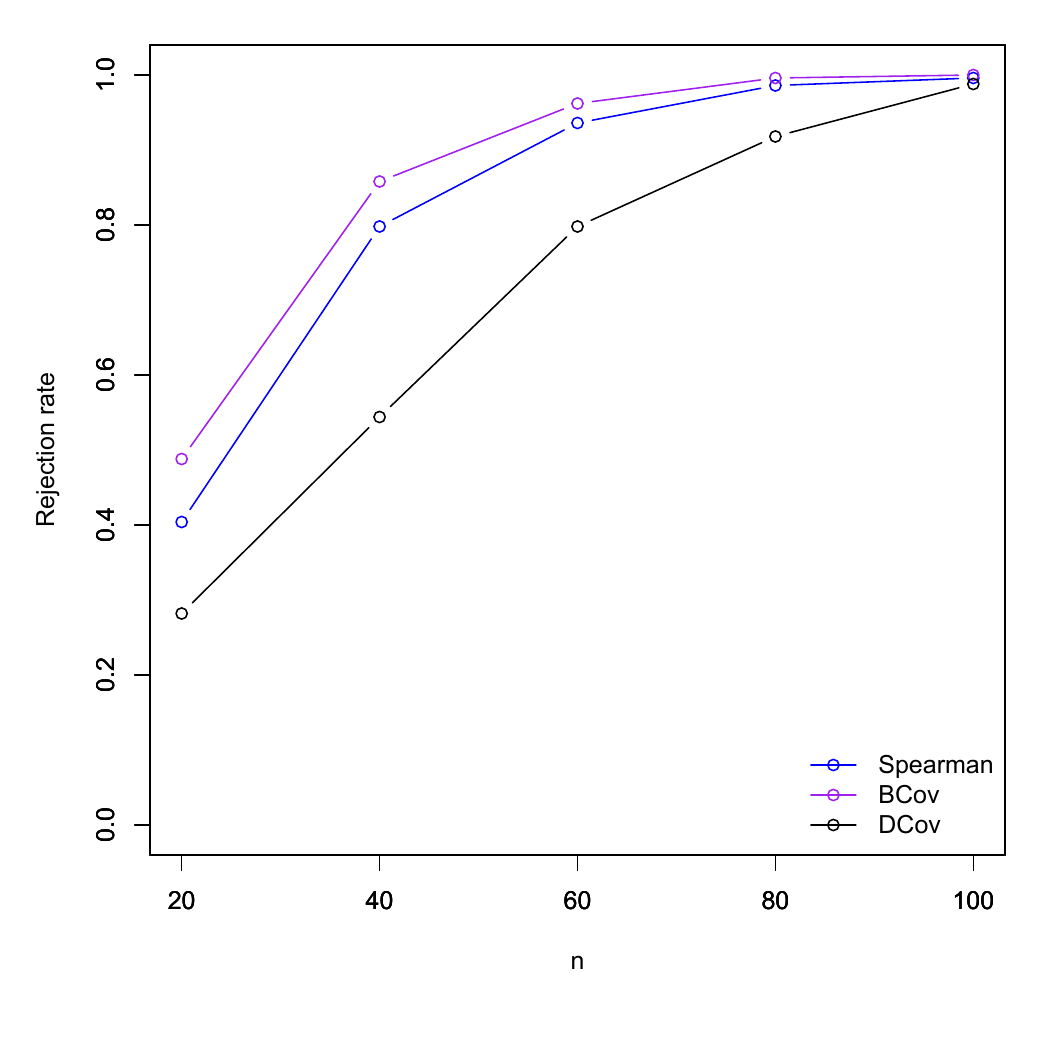}
    \end{subfigure}
    \vspace*{-30pt}
    \caption{The rejection rates of the Spearman metric rank-based, BCov and DCov tests for the spherical normal distribution (left panel) and Cauchy distribution (right panel) of $\bepsilon_i$ for different values of the sample size $n$. The red line represents the nominal level.}
  \label{plot.nDiff}
\end{figure}

\section{Conclusion and discussion}\label{Sec.Conclusion}

In Euclidean space, quantiles, ranks, and signs are crucial to nonparametric reasoning. They serve as the foundation for numerous well-known procedures for testing hypotheses. These test statistics are, hence, conceptually straightforward, computationally efficient, and theoretically beautiful. We introduced local and global metric quantiles, metric ranks, and metric signs in metric space as their isometry-invariant counterparts. In addition to exhibiting the desired theoretical qualities, their empirical counterparts are computationally efficient.  They facilitate the application of conventional nonparametric statistical inference methods to metric spaces. Consequently, these cutting-edge concepts open new research opportunities for metric space-valued data. For instance, in metric spaces, we can create a distribution-free test of independence using metric ranks or metric signs. Quantile regression, in which the covariates and/or the response are measured in metric spaces, is an additional challenging but fascinating research topic.



\bibliographystyle{Chicago}
\bibliography{MetricQF}


\end{document}